\documentclass{amsart2020}

\newtheorem{theorem}{Theorem}[section]
\newtheorem{lemma}[theorem]{Lemma}
\newtheorem{proposition}[theorem]{Proposition}

\newtheorem{problem}[theorem]{Problem}
\newtheorem{conjecture}[theorem]{Conjecture}

\theoremstyle{definition}
\newtheorem{definition}[theorem]{Definition}
\newtheorem{example}[theorem]{Example}
\newtheorem{remark}[theorem]{Remark}

\usepackage{amscd,amssymb}
\begin{document}

%\title[Weak polynomial identities of third degree \ldots and more]
%{Weak polynomial identities of third degree \ldots and more}
\title[Weak polynomial identities and their applications]
{Weak polynomial identities\\
and their applications}

\author[Vesselin Drensky]
{Vesselin Drensky}

\thanks{Partially supported by
Grant KP-06 N 32/1 of 07.12.2019 ``Groups and Rings -- Theory and Applications'' of the Bulgarian National Science Fund.}


\maketitle

\begin{center}
Institute of Mathematics and Informatics\\
Bulgarian Academy of Sciences\\
Acad. G. Bonchev Str., Block 8\\
1113 Sofia, Bulgaria\\
drensky@math.bas.bg
\end{center}

\begin{abstract}
Let $R$ be an associative algebra over a field $K$ generated by a vector subspace $V$.
The polynomial $f(x_1,\ldots,x_n)$ of the free associative algebra $K\langle x_1,x_2,\ldots\rangle$
is a weak polynomial identity for the pair $(R,V)$ if it vanishes in $R$ when evaluated on $V$.
We survey results on weak polynomial identities
and on their applications to polynomial identities and central polynomials of associative
and close to them nonassociative algebras and on the finite basis problem.
We also present results on weak polynomial identities of degree three.
\end{abstract}

\section{Introduction}

In what follows $K$ will be an arbitrary field and all algebras and vector spaces will be over $K$.
One of our main objects will be the free associative algebra $K\langle X\rangle=K\langle x_1,x_2,\ldots\rangle$ (with or without 1).
Sometimes we shall use also other symbols for its generators, e.g. $x,y,y_i$, etc.
We shall fix the characteristic of $K$ when necessary.
Let $R$ be an associative algebra (with or without 1)
and let $V$ be a vector subspace of $R$ which generates $R$ as an algebra.
The element $f(x_1,\ldots,x_n)$ of $K\langle X\rangle$
is a {\it weak polynomial identity for the pair} $(R,V)$ if
\[
f(v_1,\ldots,v_n)=0\text{ for all }v_1,\ldots,v_n\in V.
\]
The set $T(R,V)$ of all weak polynomial identities of the pair $(R,V)$ is called a {\it weak T-ideal}
and is an ideal of $K\langle X\rangle$ closed under all linear endomorphisms of $K\langle X\rangle$.
(This means that if $f(x_1,\ldots,x_n)\in T(R,V)$ and $u_1,\ldots,u_n$ are linear combinations of
elements in $X$, then $f(u_1,\ldots,u_n)\in T(R,V)$.)

Weak polynomial identities were introduced by Razmyslov \cite{Ra1, Ra2} as a method to attack two classical problems:
\begin{itemize}
\item {\it To find the polynomial identities of the (associative) algebra $M_d(K)$ of the $d\times d$ matrices over $K$.}
\item {\it To prove that for any $d>2$ there exist central polynomials for $M_d(K)$.}
\end{itemize}
See also the book by Razmyslov \cite{Ra3} for the account of his methods and results.

It is a natural question to describe the weak polynomial identities of a given pair $(R,V)$
and the generating set of the weak T-ideal $T(R,V)$.
Initially the vector subspace $V$ generating the algebra $R$ was assumed to be a Lie subalgebra
of the adjoint Lie algebra $R^{(-)}$ of $R$
with respect to the multiplication defined by the Lie bracket $[r_1,r_2]=r_1r_2-r_2r_1$, $r_1,r_2\in R$
(and $R$ is an {\it enveloping algebra} of the Lie algebra $V$).
Then the weak T-ideal $T(R,V)$ is closed under substitutions of elements of the free Lie algebra $L(X)$
canonically embedded in the free associative algebra $K\langle X\rangle$.
Later one considered also the case when $V$ is a subalgebra of of the Jordan algebra $R^{(+)}$ of $R$ with respect to the Jordan multiplication
$\displaystyle r_1\circ r_2=\frac{1}{2}(r_1r_2+r_2r_1)$, $r_1,r_2\in R$.
Finally, one may consider the case when the vector subspace $V$ does not have
any additional structure. In all these three cases one can define {\it varieties of pairs}: {\it varieties of Lie} and {\it Jordan pairs} and
{\it varieties of arbitrary pairs} or {\it L-varieties}.

The first purpose of the present paper is to survey results on weak polynomial identities
related with their applications to polynomial identities and central polynomials of associative
and close to them nonassociative algebras and on the finite basis problem.
These are the main topics on Sections \ref{section matrix identities}, \ref{section central polynomials},
\ref{section Jordan algebras} and \ref{section weak central polynomials}.
In these sections we include also related results obtained without using weak identities.
In Section \ref{section nonassociative algebras} we discuss what happens if we consider arbitrary nonassociative algebras.
In Section \ref{section L-varieties}
we survey results on L-varieties with special attention
on the case when they satisfy a weak polynomial identity of degree 3.

A significant part of the paper deals with associative, Lie and Jordan algebras satisfying polynomial identities.
For a background and further reading we refer e.g. to Drensky \cite{D8}, Drensky and Formanek \cite{DF}, Giambruno and Zaicev \cite{GZ},
Kanel-Belov, Karasik and Rowen \cite{KBKR} for associative algebras,
Bahturin \cite{Ba} for Lie algebras and Zhevlakov, Slinko, Shestakov and Shirshov \cite{ZSSS} for Jordan algebras.

Besides varieties of pairs one may consider also varieties of representations of Lie algebras and groups.
This is a topic intensively developed by the school of Plotkin in Riga. For details see e.g. \cite{Pl, PlVov} and \cite{Vov}.

\section{Weak polynomial identities and polynomial identities of matrices}\label{section matrix identities}

One of the most attracting and most important problems in the theory of algebras with polynomial identities (or PI-algebras)
is to find a generating system (or a {\it basis}) of the polynomial identities of the $d\times d$ matrix algebra $M_d(K)$.
In other words, we consider the {\it T-ideal} $T(M_d(K))$ of the polynomial identities of $M_d(K)$
(in the language of weak polynomial identities $T(M_d(K))$ coincides with the weak T-ideal $T(M_d(K),M_d(K))$ of the pair $(M_d(K),M_d(K))$).
We want to find a system of polynomial identities $f_i(x_1,\ldots,x_{n_i})\in T(M_d(K))$ which generate $T(M_d(K))$ as a T-ideal, i.e.
$T(M_d(K))$ is the minimal ideal of $K\langle X\rangle$ which contains all $f_i(x_1,\ldots,x_{n_i})$ and their substitutions
$f_i(u_1,\ldots,u_{n_i})$, $u_1,\ldots,u_{n_i}\in K\langle X\rangle$.

We shall survey some of the results for the explicit form of the polynomial identities of $M_2(K)$
and partial results on the identities of $M_d(K)$ for any $d\geq 2$. The most important case is when the field $K$ is of characteristic 0.
Razmyslov \cite{Ra1} solved this problem for $2\times 2$ matrices.
He considered the weak polynomial identities of the pair $(M_2(K),sl_2(K))$, where $sl_2(K)$ is the Lie algebra of the $2\times 2$ traceless matrices.
A typical example of a weak polynomial identity for this pair comes from the Cayley-Hamilton theorem.
If $a\in M_2(K)$, then
\[
a^2-\text{tr}(a)a+\det(a)e=0,\quad e=e_2=\left(\begin{matrix}1&0\\
0&1\\
\end{matrix}\right).
\]
If $a\in sl_2(K)$, then $\text{tr}(a)=0$ and $a^2=-\det(a)e$.
Hence $a^2$ is a scalar matrix and
\[
[a^2,b]=0\text{ for all }b\in M_2(K).
\]
Hence
\begin{equation}\label{WPI for matrices of order 2}
[x^2,y]
\end{equation}
is a weak polynomial identity for the pair $(M_2(K),sl_2(K))$.
The main results in \cite{Ra1} are the following.

\begin{theorem}\label{identities in 2x2 matrices}
Let $\text{\rm char}(K)=0$.

{\rm (i)} The weak polynomial identities of the pair $(M_2(K),sl_2(K))$ follow from the weak identity
$[x^2,y]$ (allowing substitutions in the variables by Lie elements).

{\rm (ii)} The polynomial identities of the Lie algebra $sl_2(K)$ follow from its identities of fifth degree.

{\rm (iii)} The polynomial identities of $M_2(K)$ follow from nine identities of degree $\leq 6$.

\noindent In both cases {\rm (ii)} and {\rm (iii)} the identities are explicitly given.
\end{theorem}

Depending on the field $K$ the bases of the polynomial identities of the algebras $M_2(K)$ and $sl_2(K)$
are known for fields $K$ of any characteristic with one exception.
The results are due to the work of many people starting from Wagner \cite{Wa} and continuing nowadays.

In 1936 Wagner \cite{Wa} showed that the algebra $M_2(K)$ satisfies the polynomial identity
\begin{equation}\label{Wagner-Hall identity}
[[x_1,x_2]^2,x_3].
\end{equation}
The proof of (\ref{Wagner-Hall identity}) is similar to the proof of (\ref{WPI for matrices of order 2})
because the trace of the commutator of two matrices is equal to 0.
In 1943 Hall \cite{H} showed that if a noncommutative division algebra satisfies the polynomial identity
(\ref{Wagner-Hall identity}) then it is generalized quaternion, i.e. four-dimensional over its center.
(Since the algebra $M_2({\mathbb R})$ and the quaternion algebra have the same polynomial identities
for a long period the identity (\ref{Wagner-Hall identity}) was called the Hall identity. We shall use the name {\it Wagner-Hall identity}.)
Wagner \cite{Wa} showed also an explicit polynomial identity for the algebra $M_d(K)$ of $d\times d$ matrices for any $d$.
Later his idea was further developed to produce polynomial identities for $M_d(K)$ which are simpler than the original ones
and are known as {\it identities of algebraicity}, see Theorem \ref{identities of algebraicity} below.
Recall that the {\it standard polynomial} of degree $d$ and the {\it Capelli polynomial} in $d$ alternating variables are, respectively,
\[
s_d(x_1,\ldots,x_d)=\sum_{\sigma\in S_d}\text{\rm sign}(\sigma)x_{\sigma(1)}\cdots x_{\sigma(d)},
\]
\[
c_d(x_1,\ldots,x_d;y_1,\ldots,y_{d-1})=\sum_{\sigma\in S_d}\text{\rm sign}(\sigma)x_{\sigma(1)}y_1x_{\sigma(2)}y_2\cdots y_{d-1}x_{\sigma(d)},
\]
where $S_d$ is the symmetric group acting on the symbols $1,2,\ldots,d$.
Both polynomials vanish when $x_1,\ldots,x_d$ are replaced by linearly dependent elements $r_1,\ldots,r_d$ in an algebra $R$.

\begin{theorem}\label{identities of algebraicity}
The algebra $M_d(K)$ satisfies the polynomial identities
\begin{equation}\label{first identity of algebraicity}
s_d([x,y],[x^2,y],\ldots,[x^d,y]),
\end{equation}
\begin{equation}\label{second identity of algebraicity}
c_{d+1}(1,x,x^2,\ldots,x^d;y_1,\ldots,y_d)=
\sum_{\sigma\in S_{d+1}}\text{\rm sign}(\sigma)x^{\sigma(0)}y_1x^{\sigma(1)}\cdots y_dx^{\sigma(d)},
\end{equation}
where $S_{d+1}$ acts on the symbols $0,1,2,\ldots,d$.
\end{theorem}

Both identities follow from the Cayley-Hamilton theorem. If $a\in M_d(K)$, then
it satisfies its characteristic equation
\begin{equation}\label{characteristic equation}
f_a(a)=a^d+\alpha_1a^{d-1}+\cdots+\alpha_{d-1}a+\alpha_de=0
\end{equation}
for suitable $\alpha_1,\ldots,\alpha_{d-1},\alpha_d\in K$. Hence the matrices
$e=e_d,a,a^2,\ldots,a^d$ are linearly dependent which proves the identity
(\ref{second identity of algebraicity}).
If we take the commutator of the characteristic equation (\ref{characteristic equation}) of $a$ with an arbitrary $b\in M_d(K)$,
we shall obtain
\[
[f_a(a),b]=[a^d,b]+\alpha_1[a^{d-1},b]+\cdots+\alpha_{d-1}[a,b]=0,
\]
i.e. the commutators $[a^d,b],[a^{d-1},b],\ldots,[a,b]$ are linearly dependent,
which gives the identity (\ref{first identity of algebraicity}).

In 1950 Amitsur and Levitzki  \cite{ALe} proved their famous theorem.

\begin{theorem}\label{Amitsur-Levitzki theorem}
The matrix algebra $M_d(K)$ satisfies the standard identity
\[
s_{2d}(x_1,\ldots,x_{2d})
\]
of degree $2d$. Every identity of degree $\leq 2d$ for $M_d(K)$
is equal to the standard identity up to a multiplicative constant.
The only exceptions are the cases when $d=1,2$ and $K={\mathbb F}_2$ is the field with two elements
because $M_1({\mathbb F}_2)={\mathbb F}_2$ satisfies the identity $x^2-x$
and $M_2({\mathbb F}_2)$ satisfies a nonhomogeneous identity of degree four in three variables.
\end{theorem}

In the same 1950 another important paper by Specht \cite{Sp} appeared where he stated a problem
which was one of the main driving forces in the theory for more than 30 years.

\begin{problem}\label{Specht problem}
Is it true that the polynomial identities of any PI-algebra $R$ follow from a finite number?
\end{problem}

The Specht problem is a partial case for associative algebras for the finite basis problem in universal algebra
(stated for groups by Neumann \cite{Ne} in 1937).
But its influence was so big that the name {\it Specht problem} nowadays is used also for groups and nonassociative algebras.
In particular, a variety of algebras satisfies the {\it Specht property} if it and every of its proper subvarieties
can be defined by a finite number of identities.

Let us go back to $2\times 2$ matrices. After Razmyslov established Theorem \ref{identities in 2x2 matrices} there were several attempts
to minimize and simplify the bases of the polynomial identities of $M_2(K)$ and $sl_2(K)$, $\text{char}(K)=0$. Filippov \cite{Fi}
reduced the identities of $sl_2(K)$ found by Razmyslov to one identity:

\begin{theorem}\label{theorem of Filippov}
Over a field $K$ of characteristic $0$ the polynomial identities of the Lie algebra $sl_2(K)$ follow from the identity
\begin{equation}\label{identity of Filippov}
[[[x_2,x_3],[x_4,x_1]],x_1]+[[[x_2,x_1],[x_3,x_1]],x_4].
\end{equation}
\end{theorem}

Using a computer, in the 1970s Rosset, see \cite{Ler}, showed that the polynomial identities of degree 5 for $M_2(K)$, $\text{char}(K)=0$,
follow from the standard identity $s_4(x_1,x_2,x_3,x_4)$ and the Wagner-Hall identity (\ref{Wagner-Hall identity}).
Bui \cite{Bu} reduced the basis of Razmyslov from Theorem \ref{identities in 2x2 matrices} to $s_4(x_1,x_2,x_3,x_4)$,
(\ref{Wagner-Hall identity}) and two more identities of degree 6.
Finally, the author \cite{D4} found a minimal basis of the identities of $M_2(K)$, $\text{char}(K)=0$.
The proof is based on representation theory of the symmetric and the general linear groups.
It uses the result of Razmyslov that the polynomial identities of $M_2(K)$ follow from those of degree $\leq 6$
but does not use the explicit form of the nine identities of the basis of Razmyslov.

\begin{theorem}\label{minimal basis of identities of 2x2 matrices}
Over a field $K$ of characteristic $0$ the polynomial identities of the algebra $M_2(K)$ follow from the identities
\[
s_4(x_1,x_2,x_3,x_4)\text{ and }[[x_1,x_2]^2,x_1].
\]
\end{theorem}

In the same paper \cite{D4} the author obtained the following basis of the polynomial identities of the Lie algebra $sl_2(K)$, $\text{char}(K)=0$:
\[
\sum_{\sigma\in S_3}\text{\rm sign}(\sigma)[x_{\sigma(1)},x_1,x_1,x_{\sigma(2)},x_{\sigma(3)}]
\]
and the Lie standard identity
\[
x_1s_4(\text{ad}(x_1),\text{ad}(x_2),\text{ad}(x_3),\text{ad}(x_4))
=\sum_{\sigma\in S_3}\text{\rm sign}(\sigma)[x_1,x_{\sigma(1)},x_{\sigma(2)},x_{\sigma(3)},x_{\sigma(4)}].
\]
(The Lie brackets are left normed: $[u_1,\ldots,u_{n-1},u_n]=[[u_1,\ldots,u_{n-1}],u_n]$, $n\geq 3$,
and $u_1\text{ad}(u_2)=[u_1,u_2]$.)

Razmyslov had a conjecture that for $d>2$ the polynomial identities of $M_d(K)$, $\text{char}(K)=0$, follow from
the standard identity $s_{2d}(x_1,\ldots,x_{2d})$
and the identity of algebraicity (\ref{first identity of algebraicity}) in Theorem \ref{identities of algebraicity}.
The result of Theorem \ref{minimal basis of identities of 2x2 matrices} gives that this conjecture holds for $d=2$.
The problem for the description of the bases of the polynomial identities of $M_d(K)$, $\text{char}(K)=0$, $d>2$, is still open even for $d=3$.
Leron \cite{Ler} proved that for $d>2$ all polynomial identities of degree $2d+1$ for $M_d(K)$ follow from the standard identity
$s_{2d}(x_1,\ldots,x_{2d})$. For $d=3$ Drensky and Kasparian \cite{DK2} made one more step and showed that the same holds for
the polynomial identities of degree 8. The computations were done by hand, without using any computers, and were based on representation theory
of the general linear group. Combining theoretical results (representation theory of the general linear group again)
and computer calculations, Benanti, Demmel, Drensky and Koev \cite{BDDK}
showed that the identities of degree $2d+2$ for $M_d(K)$ follow from the standard identity $s_{2d}(x_1,\ldots,x_{2d})$ also for $d=4,5$.
Drensky and Kasparian \cite{DK1} showed that for $d\geq 3$ the identity (\ref{second identity of algebraicity})
does not follow from $s_{2d}(x_1,\ldots,x_{2d})$ and (\ref{first identity of algebraicity}).
Okhitin \cite{Ok} described all polynomial identities for $M_3(K)$ of degree 9 in two variables
and found one which does not follow from $s_6(x_1,\ldots,x_6)$ and (\ref{second identity of algebraicity}).
Domokos \cite{Do} exhibited another three new polynomial identities for $M_3(K)$ which are of degree 9 and in three variables.
See \cite{BDDK} for a survey on other computational results on polynomial identities of matrices in characteristic 0.

The picture and the methods for the polynomial identities of matrices in positive characteristic are quite different from the case of characteristic 0.
If in characteristic 0 many of the investigations are based on representation theory of the symmetric and the general linear group,
the methods in positive characteristic depend also on the number of elements of the field.
When $K$ is a finite field then one uses a lot of structure ring theory and for infinite fields combinatorial methods are combined
with characteristic free invariant theory.

In 1964 Oates and Powell \cite{OP} proved that the variety generated by a finite group is finitely based and has a finite number of subvarieties.
In 1973 the methods of \cite{OP} were transferred successfully by Kruse \cite{Kr} and L'vov \cite{Lv1} and they established the result
for finite associative rings and for finite dimensional associative algebras over finite fields.
The study of polynomial identities of (not only associative) finite rings and algebras was very popular in the theory of PI-algebras in
the 1970s and 1980s.

\begin{theorem}\label{Theorem of Kruse and L'vov}
The polynomial identities of every finite associative ring and of every finite dimensional associative algebra over a finite field
follow from a finite number of identities.The variety generated by such a ring or an algebra has a finite number of subvarieties
which also are generated by a finite object.
\end{theorem}

In the special case when $K={\mathbb F}_q$ is a finite field the polynomial identities of $M_1({\mathbb F}_q)={\mathbb F}_q$ follow from the identity
$x^q-x$. The proof uses the theorem of Jacobson, see e.g. \cite[Theorem 3.1.2]{He}, that if for any element $a$ in a ring $R$ there exists
an integer $n=n(a)>1$ such that $a^n=a$, then $R$ is commutative.
The identities of $M_2({\mathbb F}_q)$ were described by Mal'tsev and Kuz'min \cite{MaK}.

\begin{theorem}\label{Theorem of Maltsev and Kuzmin}
The following polynomial identities form a basis of the polynomial identities of $M_2({\mathbb F}_q)$
\[
(x-x^q)(y-y^q)(1-[x,y]^{q-1}),
\]
\[
(x-x^q)\circ(y-y^q)-((x-x^q)\circ(y-y^q))^q.
\]
(In the case of fields of characteristic $2$ we assume that $u_1\circ u_2=u_1u_2+u_2u_1$.)
\end{theorem}

Later Genov \cite{Ge} and Genov and Siderov \cite {GeS} found, respectively,
bases of the polynomial identities of $M_3({\mathbb F}_q)$ and $M_4({\mathbb F}_q)$.

Bakhturin and Ol'shanski\u{\i} \cite{BaO} developed further the methods of Oates and Powell \cite{OP},
Kruse \cite{Kr} and L'vov \cite{Lv1} involving also quasivarieties and proved analogues
for finite Lie rings and finite dimensional Lie algebras over finite fields. Then Semenov \cite{Se} found
a basis of the identities of $sl_2({\mathbb F}_q)$, $\text{char}({\mathbb F}_q)>3$
(and to the best of our knowledge the case $\text{char}({\mathbb F}_q)=3$ is still open):

\begin{theorem}\label{Theorem of Semenov}
The following polynomial identities form a basis of the polynomial identities of $sl_2({\mathbb F}_q)$,
$\text{\rm char}({\mathbb F}_q)>3$:
\[
[x,y](1-\text{\rm ad}^{q^2-1}(x)-\text{\rm ad}^{q-1}(y)+\text{\rm ad}^{q^2-1}(x)\text{\rm ad}^{q-1}(y))
\]
\[
+[x,y](\text{\rm ad}^{q^2}(x)-\text{\rm ad}(x))\text{\rm ad}^{q-2}([x,y])(\text{\rm ad}^{q^2}(y)-\text{\rm ad}(y))
\]
\[
-y((\text{\rm ad}^{q^2}(x)-\text{\rm ad}(x))\text{\rm ad}(x))^q(\text{\rm ad}^{q^2-2}(y)-\text{\rm ad}^{q-2}(y)),
\]
\[
[y,x,x](\text{\rm ad}^{q^2}(x)-\text{\rm ad}(x)).
\]
\end{theorem}

When $\text{char}(K)=2$ the algebra $sl_2(K)$ is nilpotent and the Lie algebra $gl_2(K)=M_2^{(-)}(K)$ all $2\times 2$ matrices satisfies
the center-by-metabelian identity
\[
[[[x_1,x_2],[x_3,x_4]],x_5].
\]
The basis of the polynomial identities of $gl_2({\mathbb F}_q)$, $q=2^m$, was found by the author \cite{D3}:

\begin{theorem}\label{gl(2) over a finite field of char 2}
The following polynomial identities form a basis of the polynomial identities of $gl_2({\mathbb F}_q)$,
$\text{\rm char}({\mathbb F}_q)=2$:
\[
[[[x_1,x_2],[x_3,x_4]],x_5],
\]
\[
[[x_1,x_2],[x_3,x_4]]+[[x_1,x_3],[x_2,x_4]]+[[x_1,x_4],[x_2,x_3]].
\]
(This identity is equal to the standard identity $s_4(x_1,x_2,x_3,x_4)$
which in characteristic $2$ can be expressed as a Lie element.)
\[
[[[x_1,x_2],x_3],[x_1,x_2]],
\]
\[
[[x_1,x_2]\text{\rm ad}^{q-1}(x_1),[x_1,x_3]]-[[x_1,x_2],[x_1,x_3]],
\]
\[
[[x_1,x_2],x_3(1+\text{\rm ad}^{q-1}(x_1))(1+\text{\rm ad}^{q-1}(x_2))(1+\text{\rm ad}^{q-1}(x_3))]-[x_1,x_2]\text{\rm ad}^q(x_3).
\]
\end{theorem}

In characteristic 2 there is another simple three-dimensional Lie algebra which is an analogue of the Lie algebra ${\mathbb R}^3$
equipped with the cross product (or the vector multiplication). It would be interesting to find a basis of the polynomial identities of this Lie algebra.

Finally, we shall survey the results for the polynomial identities of $M_2(K)$ and $sl_2(K)$ when $K$ is an infinite field of positive characteristic.
Partial results were obtained by Koshlukov \cite{Ko2} and then
Colombo and Koshlukov \cite{CoKo} obtained a minimal basis of the identities of $M_2(K)$ when $K$ is an infinite field of odd characteristic.

\begin{theorem}\label{basis of Colombo and Koshlukov}
Let $K$ be an infinite field of characteristic $p>2$. If $\text{\rm char}(K)>3$, then the polynomial identities of $M_2(K)$ follow from
the identities
\[
s_4(x_1,x_2,x_3,x_4)\text{ and }[[x_1,x_2]\circ[x_3,x_4],x_5].
\]
If $\text{\rm char}(K)=3$, then one has to add the polynomial identity
\[
2[x_1,x_2]\circ(u\circ v)-[x_1,u,v,x_2]-[x_1,v,u,x_2]+[x_2,u,x_1,v]+[x_2,v,x_1,u].
\]
\end{theorem}

The case for the polynomial identities of $M_2(K)$ over an infinite field of characteristic 2 is still open.
In virtue of the result for the polynomial identities of $M_2(K)$ considered as a Lie algebra, see Theorem \ref{Lie identities in char 2},
very probably the following conjecture is true.

\begin{conjecture}\label{conjecture for char 2}
The polynomial identities of $M_2(K)$ over an infinite field of characteristic $2$ do not follow from a finite number.
\end{conjecture}

The polynomial identities of the Lie algebra $sl_2(K)$ over an infinite field of odd characteristic were described by Vasilovskij \cite{Va1}.

\begin{theorem}\label{sl2 for infinite field of odd characteristic}
When $K$ is an infinite field of characteristic $p>2$ the polynomial idenetities of the Lie algebra $sl_2(K)$ follow from
the identity of Filippov (\ref{identity of Filippov}).
\end{theorem}

As in the case of finite fields of characteristic 2 instead of the algebra $sl_2(K)$ which is nilpotent of class 2 one considers
the algebra $gl_2(K)$. In 1970 Vaughan-Lee \cite{VL1} proved that this algebra does not have a finite basis of its polynomial identities.
In his Ph.D. Thesis from 1979 \cite{D2} the author found an explicit infinite basis of these identities
but the result was not published in a journal paper. Later the same result was obtained independently by Lopatin \cite{Lo}.

\begin{theorem}\label{Lie identities in char 2}
Let $K$ be an infinite field of characteristic $2$.
Then the polynomial identities of the Lie algebra $gl_2(K)$ do not have a finite basis and follow from the identities
\[
[[[x_1,x_2],[x_3,x_4]],x_5],
\]
\[
[[x_1,x_2,x_3,\ldots,x_n],[x_1,x_2]],\quad n=3,4,\ldots,
\]
\[
[[x_4,x_1,x_5,\ldots,x_n],[x_2,x_3]]+[[x_3,x_1,x_5,\ldots,x_n],[x_4,x_2]]
\]
\[
+[[x_2,x_1,x_5,\ldots,x_n],[x_3,x_4]],
\quad n=4,5,\ldots.
\]
\end{theorem}

Vaughan-Lee \cite{VL1} proved that a basis of the identities of $gl_2(K)$ consists of the center-by-metabelian identity, the identities
$[[x_1,x_2,x_3,\ldots,x_n],[x_1,x_2]]$
and multilinear identities which do not have consequences of higher degree modulo the center-by-metabelian identity.
The author \cite{D2} used a theorem of Kuz'min \cite{Ku} to find these multilinear identities. Lopatin \cite{Lo} described these multilinear identities
with direct computations modulo the center-by-metabelian identity.

As in the case of finite fields of characteristic 2 it is interesting to find the polynomial identities in the three-dimensional analogue of
${\mathbb R}^3$.

\begin{conjecture}
For an infinite field $K$ of characteristic $2$ the polynomial identities of the three-dimensional simple Lie algebra $K^3$
equipped with the cross product do not follow from a finite number.
\end{conjecture}

We shall conclude this section with some comments on the solution of different versions of the Specht problem for associative and Lie algebras.
For further reading we refer to the book by Kanel-Belov, Karasik and Rowen \cite{KBKR}.

In 1987 Kemer \cite{Ke1}, see his book \cite{Ke3} for an account,
solved into affirmative the original Specht problem \cite{Sp} for algebras over a field of characteristic 0.
Then he proved a weaker version of the Specht problem for finitely generated algebras over an infinite field \cite{Ke2}.

\begin{theorem}\label{results of Kemer}
{\rm (i)} Every T-ideal of the free associative algebra $K\langle X\rangle$
over a field of characteristic $0$ has a finite basis.

{\rm (ii)} The polynomial identities in $d$ variables in every T-ideal of the finitely generated free associative algebra $K\langle x_1,\ldots,x_d\rangle$
over an infinite field follow from a finite number.
\end{theorem}

Later Theorem \ref{results of Kemer} (ii) was generalized by Belov \cite{B3}
for T-ideals of finitely generated associative algebras over commutative associative noetherian rings,
see also Belov-Kanel, Rowen and Vishne \cite{BKRV} for detailed exposition.

On the other hand, the Specht problem has a negative solution over fields of positive characteristic.
The first counterexamples were given in 1999 by Belov \cite{B1, B2}.

\begin{theorem}\label{counterexample of Belov}
Over an arbitrary field of positive characteristic there are T-ideals which do not have finite bases of their polynomial identities.
\end{theorem}

In the same 1999 Grishin \cite{Gr1, Gr2} and in 2002 Gupta and Krasilnikov \cite{GuKr}
constructed simple counterexamples to the Specht problem in characteristic 2.

\begin{theorem}\label{Grishin, Gupta-Krasilnikov}
Let $K$ be a field of characteristic $2$. Then the following systems of polynomial identities do not follow from a finite number:

{\rm (i)} The system of Grishin:
\[
y_1^4z_1^4x_1^2\cdots x_n^2z_2^4y_2^4y_1^4z_1^4x_{n+1}^2\cdots x_{2n}^2z_2^4y_2^4,\quad n=0,1,2,\ldots;
\]

{\rm (ii)} The system of Gupta and Krasilnikov:
\[
[x,y^2]x_1^2\cdots x_n^2[x,y^2]^3,\quad n=0,1,2,\ldots .
\]
\end{theorem}

The example of Grishin in Theorem \ref{Grishin, Gupta-Krasilnikov} gives that the algebras satisfying his polynomial identities
satisfy the polynomial identity $x^{32}$ and hence are nil of bounded index. By the theorem of Levitzki \cite{Le}
such algebras are locally nilpotent (i.e. finitely generated algebras are nilpotent).

The Specht problem is still open for Lie algebras over a field of characteristic 0.
In the positive direction we shall state the following result of Il'tyakov \cite{I2}, see also \cite{I3} for a more detailed exposition.

\begin{theorem}\label{Theorem of Iltyakov for Lie algebras}
{\rm (i)} If $L$ is a finitely generated Lie algebra over a field of characteristic $0$ and if $L$ has an enveloping algebra $R$ which is a PI-algebra,
then every pair satisfying all weak polynomial identities of the pair $(R,L)$ has a finite basis of its weak identities.

{\rm (ii)} If $L$ is a finitely generated Lie algebra and its algebra of multiplications $\text{\rm ad}(L)$ is an (associative) PI-algebra,
then the variety of Lie algebras generated by $L$ satisfies the Specht property.
In particular, this holds when the Lie algebra $L$ is finite dimensional.
\end{theorem}

\section{Weak polynomial identities and central polynomials}\label{section central polynomials}

\begin{definition}\label{definition of central polynomials}
The polynomial $c(x_1,\ldots,x_n)\in K\langle X\rangle$ is a {\it central polynomial} for the associative algebra $R$,
if $c(x_1,\ldots,x_n)$ is not a polynomial identity for $R$ and for all $r_1,\ldots,r_n\in R$ it holds that
$c(r_1,\ldots,r_n)$ belongs to the center of $R$.
\end{definition}

An example of a central polynomial for the algebra $M_2(K)$ of $2\times 2$ matrices is
\[
c(x_1,x_2)=[x_1,x_2]^2.
\]
For the proof the Wagner-Hall identity (\ref{Wagner-Hall identity}) gives that
$c(x_1,x_2)$ has only central values when evaluated on $M_2(K)$ and it is easy to see that it is not a polynomial identity
because there are commutators $[a_1,a_2]$, $a_1,a_2\in M_2(K)$ which have nonzero determinants.

In a talk given in 1956 Kaplansky \cite{K1} (see also the revised version \cite{K2} from 1970) asked several problems which
motivated significant research activity. One of the problems is for the existence of central polynomials for
the $d\times d$ matrix algebra $M_d(K)$, $d>2$.

Latyshev and Shmelkin \cite{LS} constructed a central (clearly non-homogeneous) polynomial in one variable
for the matrix algebra $M_d({\mathbb F}_q)$ over a finite field ${\mathbb F}_q$.
The first central polynomials over an arbitrary ground field were constructed by Formanek \cite{F1} and Razmyslov \cite{Ra2}
and this gave rise to a serious revision of the theory of algebras with polynomial identities, see e.g. the books
by Procesi \cite{P}, Jacobson \cite{J}, Rowen \cite{Ro}, Formanek \cite{F4}, Drensky and Formanek \cite{DF}, Giambruno and Zaicev \cite{GZ}.
In particular, important theorems on PI-algebras were established or simplified using central polynomials.
Later Kharchenko \cite{Kh} gave a short proof for the existence of central polynomials for $M_d(K)$ using classical results of Amitsur \cite{A1, A2}.

The following theorem of Formanek \cite{F1} gives the existence of central polynomials for $M_d(K)$ over an arbitrary field $K$.

\begin{theorem}\label{central polynomial of Formanek}
For any $d\geq 1$ and any field $K$ the algebra $M_d(K)$ has a central polynomial of degree $d^2$.
\end{theorem}

The central polynomials of Formanek are obtained with the following construction.
Let $K[u_1,\ldots,u_{d+1}]$ be the polynomial algebra in $d+1$ commuting variables.
We define a linear mapping $\theta$ (not an algebra homomorphism)
from $K[u_1,\ldots,u_{d+1}]$ to the free algebra $K\langle x,y_1,\ldots,y_d\rangle$
of rank $d+1$ in the following way. If
\[
g(u_1,\ldots,u_{d+1})=\sum\alpha_nu_1^{n_1}\cdots u_{d+1}^{n_{d+1}},
\quad\alpha_n\in K,
\]
then
\[
\theta(g)(x,y_1,\ldots,y_d)=
\sum\alpha_ax^{n_1}y_1x^{n_2}y_2x^{n_3}y_3
\cdots x^{n_d}y_dx^{n_{d+1}}.
\]
Let
\[
g(u_1,\ldots,u_{d+1})=
\prod_{2\leq i\leq d}(u_1-u_i)(u_{d+1}-u_i)
\prod_{2\leq i<j\leq d}(u_i-u_j)^2.
\]
Then the polynomial of the free associative algebra
$K\langle x,y_1,\ldots,y_d\rangle$
\[
c(x,y_1,\ldots,y_d)=\theta(g)(x,y_1,\ldots,y_d)+
\theta(g)(x,y_2,\ldots,y_d,y_1)
\]
\[
+\cdots+\theta(g)(x,y_d,y_1,\ldots,y_{d-1})
\]
is the central polynomial for the matrix algebra $M_d(K)$ from Theorem \ref{central polynomial of Formanek}.

The approach of Razmyslov \cite{Ra2} is based on weak polynomial identities. The key role plays the following map
called the {\it Razmyslov transform}.
Let the polynomial
$f(x,y_1,\ldots,y_n)\in K\langle x,y_1,\ldots,y_n\rangle$
be linear (i.e. homogeneous of degree 1) in the variable $x$.
We write $f$ in the form
\[
f=\sum g_ixh_i,\quad g_i,h_i\in K\langle y_1,\ldots,y_n\rangle,
\]
and define the Razmyslov transform of $f$ as
\[
f^{\ast}(x,y_1,\ldots,y_n)=\sum h_ixg_i.
\]

For example, if
\[
f(x,y_1,y_2)=[xy_1+y_1x,y_2]=
1\cdot x\cdot y_1y_2+y_1\cdot x\cdot y_2
-y_2\cdot x\cdot y_1-y_2y_1\cdot x\cdot 1,
\]
then
\[
f^{\ast}(x,y_1,y_2)=
y_1y_2\cdot x\cdot 1+y_2\cdot x\cdot y_1
-y_1\cdot x\cdot y_2-1\cdot x\cdot y_2y_1
\]
\[
=[y_2,x]y_1+y_1[y_2,x]=2[y_2,x]\circ y_1.
\]
Now the central polynomials are obtained using the following lemma.

\begin{lemma}\label{lemma of Razmyslov}
Let the polynomial
\[
f(x,y_1,\ldots,y_n)\in K\langle x,y_1,\ldots,y_n\rangle
\]
be homogeneous of first degree in $x$ and let
$f^{\ast}(x,y_1,\ldots,y_n)$ be the polynomial obtained
by the Razmyslov transform. Then:

{\rm (i)} $f(x,y_1,\ldots,y_n)$ is a polynomial identity for
the matrix algebra $M_d(K)$
if and only if  $f^{\ast}(x,y_1,\ldots,y_n)$ is a polynomial
identity.

{\rm (ii)} $f^{\ast}(x,y_1,\ldots,y_n)$ is a central polynomial for
the algebra $M_d(K)$
if and only if $f(x,y_1,\ldots,y_n)$ is a weak polynomial identity for the pair $(M_d(K),sl_d(K))$ and
$f([x,y_0],y_1,\ldots,y_n)$
is a polynomial identity for $M_d(K)$.
\end{lemma}

\begin{example}\label{central polynomial for M2}
If we linearize the weak polynomial identity $[x^2,y]$ for the pair $(M_2(K),sl_2(K))$,
we shall obtain the weak polynomial identity $[x\circ y_1,y]$. Then we replace $y_1$ by $[y_1,y_2]$ and obtain the polynomial
\[
f(x,y_1,y_2,y_3)=[x\circ [y_1,y_2],y_3].
\]
Now $f([x,y_0],y_1,y_2,y_3)=[[x,y_0]\circ [y_1,y_2],y_3]$ is a polynomial identity for $M_2(K)$ and hence
\[
f^{\ast}(x,y_1,y_2,y_3)=[x\circ [y_1,y_2],y_3]=[y_3,x]\circ [y_1,y_2]
\]
is a central polynomial for $M_2(K)$.
\end{example}

The central polynomials of Razmyslov \cite{Ra2} were obtained using the Capelli polynomial.

\begin{theorem}\label{central polynomial of Razmyslov}
Let
\[
f=f(x,z_1,\ldots,z_{2d^2-2},y_1,\ldots,y_{d^2-1})
\]
\[
=c_{d^2}(x,[z_1,z_2],\ldots,[z_{2d^2-3},z_{2d^2-2}];
y_1,\ldots,y_{d^2-1})
\]
where
$c_n(x_1,\ldots,x_n;y_1,\ldots,y_{n-1})$ is the Capelli
polynomial. The Razmyslov transform applied to $f$
gives a multilinear central polynomial for $M_d(K)$ of degree $3d^2-2$ over
any field $K$.
\end{theorem}

Halpin \cite{Ha} used the Cayley-Hamilton theorem to construct another weak polynomial identity for the pair $(M_d(K),sl_d(K))$.
The matrix $a\in M_d(K)$ satisfies its characteristic equation
\[
a^d+\alpha_1a^{d-1}+\cdots+\alpha_{d-1}a_1+\alpha_de=0
\]
and $\alpha_1=-\text{tr}(a)$. If $a\in sl_d(K)$, then $\text{tr}(a)=0$
and $a^d,a^{d-2},a^{d-3},\ldots,a,e$ are linearly dependent. Hence the Capelli identity
\[
f(x,y_1,\ldots,y_{d-1})=c_d(1,x,x^2,\ldots,x^{d-2},x^d;y_1,\ldots,y_{d-1})
\]
is a weak polynomial identity. It is of degree $\displaystyle \frac{1}{2}(d^2-d+2)$ in $x$
and it is easy to see that it is not a polynomial identity.
The complete linearization of $f(x,y_1,\ldots,y_{d-1})$ in $x$ gives a polynomial
$g(x,z_1,\ldots,z_k;y_1,\ldots,y_{d-1})$, $\displaystyle k=\frac{1}{2}(d-1)d$. As in Example \ref{central polynomial for M2},
$g(x,[z_1,z_2],\ldots,[z_{2k-1},z_{2k}];y_1,\ldots,y_{d-1})$ is a weak polynomial identity of degree $d^2$ which becomes a polynomial identity
when $x$ is replaced by a commutator. Applying Lemma \ref{lemma of Razmyslov} we obtain a central polynomial of degree $d^2$.

Hence both methods of Formanek and Razmyslov produce central polynomials for $M_d(K)$ of degree $d^2$.
It is easy to see that the minimal degree of the central polynomials for $M_1(K)=K$ and $M_2(K)$ is, respectively, equal to 1 and 4.
There was a conjecture that {\it the minimal degree of the central polynomials for $M_d(K)$ is $d^2$.}
For $d=3$ Drensky and Azniv Kasparian \cite{DK3} proved that the minimal degree
of the central polynomials is 8 (with computations by hand, without computers). Formanek \cite{F4} stated the following conjecture.

\begin{conjecture}
The minimal degree of the central polynomials for $M_d(K)$ over a field $K$ of characteristic $0$ is equal to
$\displaystyle \frac{1}{2}(d^2+3d-2)$.
\end{conjecture}

Drensky and Piacentini-Cattaneo \cite{DPC} constructed a central polynomial of degree 13
for the algebra $M_4(K)$ (but we still do not know whether $M_4(K)$ has central polynomials of lower degree).
Later Drensky \cite{D7} extended this construction to central polynomials of degree $(d-1)^2+4$ for all $d>2$
(but $(d-1)^2+4>(d^2+3d-2)/2$ for $d>4$).

In all examples above the degree of the weak polynomial identity is lower than the degree of the central polynomial
produced by the method of Razmyslov and this is a big advantage from computational point of view.
Starting with the weak polynomial identity of degree 3 in Example \ref{central polynomial for M2} we obtain a central polynomial of degree 4.
Drensky and Rashkova \cite{DR} described the weak polynomial identities of degree 6 for the pair $(M_3(K),sl_3(K))$
and one of them produced a central polynomial of degree 8, as in \cite{DK3}.
The central polynomial of degree 13 for $M_4(K)$ in \cite{DPC} was obtained from a weak polynomial identity of degree 9.
The results in \cite{DR} and \cite{DPC} were obtained as a combination of computer search and representation theory of the general linear group.
For the central polynomial of degree $(d-1)^2+4$ for $M_d(K)$, $d>2$ in \cite{D7}
the author used a weak polynomial identity of degree $\displaystyle \frac{1}{2}(d^2-d+6)$ which was similar to those in \cite{DR} and \cite{DPC}.

Other central polynomials with different additional properties were constructed by several authors.
Confirming a conjecture of Regev \cite{Re} Formanek \cite{F3} showed that a polynomial in two sets of $d^2$ skew-symmetric variables
is a central polynomial for $M_d(F)$. The existence of such a polynomial has important consequences for the study
of the sequence of cocharacters of matrices, see \cite{F2}. Giambruno and Valenti \cite{GV} described other central polynomials with a similar property.

\section{Weak polynomial identities and Jordan algebras}\label{section Jordan algebras}

Given an associative algebra $R$ over a field $K$ of characteristic different from 2,
we shall equip it with the structure of Jordan algebra (denoted by $R^{(+)}$) with respect to the Jordan multiplication
$\displaystyle r_1\circ r_2=\frac{1}{2}(r_1r_2+r_2r_1)$.
Considering the applications of weak polynomial identities to Jordan algebras, most of the investigations concern two algebras.
One of them is the {\it Jordan algebra $H_d(K)$ of the symmetric $d\times d$ matrices}.
The other algebra is the {\it Jordan algebra of a nondegenerate symmetric bilinear form} of the $d$-dimensional vector space $V_d$,
$d=2,3,\ldots, \infty$. When the field $K$ is algebraically closed all nondegenerate symmetric bilinear forms are equivalent.
In the general case $V_d$ may have several nonequivalent nondegenerate symmetric bilinear forms.
This implies that there are several nonisomorphic Jordan algebras
associated with the forms. When the field $K$ is infinite, all these algebras have the same polynomial identities.
We shall fix one of the forms of $V_d$ and shall denote the corresponding Jordan algebra by $B_d$.
As a vector space $B_d$ is the direct sum of $K$ and $V_d$.
The multiplication of $B_d$ is defined by
\[
(\alpha_1+a_1)(\alpha_2+a_2)=(\alpha_1\alpha_2+\langle a_1,a_2\rangle)+(\alpha_1a_2+\alpha_2a_1),
\]
where $\alpha_1,\alpha_2\in K$, $a_1,a_2\in V_d$ and $\langle a_1,a_2\rangle$ is the value of the nondegenerate symmetric bilinear form of $V_d$.
The algebra $B_d$ is a subalgebra of the {\it Clifford algebra} $C_d=C(V_d)$ on the vector space $V_d$ considered as a Jordan algebra.
Recall that if the vector space $V_d$ has a basis $\{v_1,\ldots,v_d\}$, then $C_d$ is the unitary associative algebra generated by $\{v_1,\ldots,v_d\}$
with defining relations
\[
v_i\circ v_j=\langle v_i,v_j\rangle,\quad 1\leq i,j\leq d.
\]
As a vector space $C_d$ has a basis consisting of all products
\[
v_{i_1}\cdots v_{i_n},\quad 1\leq i_1<\cdots<i_n\leq d.
\]

Below we shall survey some results on the bases of the weak and ordinary polynomial identities of $H_2(K)$ and $B_d$.
For more details and additional results we refer to the corresponding papers and the references there.

The Dniester Notebook \cite{Dn} is one of the main sources of open problems in the theory of polynomial identities of Jordan algebras.
In \cite[Problem 2.96]{Dn} Slinko asked the following problem:

\begin{problem}\label{problem of Slinko}
Find a basis of the weak identities of the pair $(M_2(K),H_2(K))$.
Do they all follow from the standard identity $s_4(x_1,x_2,x_3,x_4)$?
\end{problem}

The answer was given by the author in \cite{D5}.

\begin{theorem}\label{WPI for symmetric matrices of order 2}
When $K$ is a field of characteristic zero all weak polynomial identities of the pair $(M_2(K),H_2(K))$
follow (as Jordan consequences) from the standard identity $s_4(x_1,x_2,x_3,x_4)$
and the metabelian identity $[[x_1,x_2],[x_3,x_4]]$.
\end{theorem}

Several problems were stated by Shestakov \cite[Problems 2.126 and 2.127]{Dn}.

\begin{problem}\label{problems of Shestakov}
{\rm (i)} Find a basis of identities of the Jordan algebra $B_d$, $d=2,3,\ldots,\infty$,
of a bilinear form over an infinite field.
Does this algebra generate a Specht variety?

{\rm (ii)} Are the varieties generated by the Jordan algebras $H_d(K)$
and $M_d^{(+)}(K)$ finitely based or Specht?
\end{problem}

The answer of the first part of Problem \ref{problems of Shestakov} (i) was given by Vasilovskij \cite{Va2}.

\begin{theorem}\label{basis of Vasilovskij}
Let $K$ be an infinite field of characteristic different from $2$. If $\text{\rm char}(K)\not=3,5,7$,
then the polynomial identities of the algebra $B_{\infty}$ follow from the identities
\[
([x,y]^2,z,t)\text{ and }
\sum_{\sigma\in S_3}\text{\rm sign}(\sigma)(x_{\sigma(1)},(x_{\sigma(2)},x,x_{\sigma(3)}),x).
\]
For the basis of the identities of the algebra $B_d$, $d<\infty$, one has to add the identities
\[
\sum_{\sigma\in S_{d+1}}\text{\rm sign}(\sigma)(x_{\sigma(1)},y_1,x_{\sigma(2)},\ldots,y_{d},x_{\sigma(d+1)}),
\]
\[
\sum_{\sigma\in S_{d+1}}\text{\rm sign}(\sigma)(x_{\sigma(1)},y_1,x_{\sigma(2)},\ldots,y_{d-1},x_{\sigma(d)})(y_d,x_{\sigma(d+1)},y_{d+1}).
\]
When $\text{\rm char}(K)=3,5,7$ additional explicit identities are given to complete the basis.
\end{theorem}

Here the multiplications are in the free Jordan algebra and
\[
(x,y,z)=(x,y)z-x(yz)
\]
is the {\it associator}.
The square of the commutator $[x,y]^2$ can be expressed as a Jordan element because
\[
[x,y]^2=4((x\circ x)\circ (y\circ y)+(x\circ y)\circ(x\circ y)-((x\circ x)\circ y)\circ y-((y\circ y)\circ x)\circ x)
\]
in the free associative algebra.

The second part of Problem \ref{problems of Shestakov} (i) was solved into affirmative for unitary Jordan algebras over a field of characteristic 0.
Il'tyakov \cite{I1} proved that the subvarieties of the variety of Jordan algebras $\text{var}(B_d)$ generated by $B_d$, $d<\infty$,
and $\text{var}(B_d)$ itself have finite bases of their identities.
Koshlukov \cite{Ko1} showed that the T-ideals containing the T-ideal of the identities of $B_{\infty}$ satisfy the ascending chain condition
which in combination with the basis of $B_{\infty}$ found by Vasilovskij \cite{Va2} establishes the Specht property of $T(B_{\infty})$.
The result of Koshlukov uses the description of the relatively free algebra $F(\text{var}(B_{\infty}))=J(X)/T(B_{\infty})$
given by the author of the present paper in \cite{D6} in the language of representation theory of the symmetric and the general linear groups.

\begin{theorem}\label{Specht property for Jordan algebra of bilinear form}
Let $K$ be a field of characteristic zero. Then every T-ideal of the free unitary Jordan algebra $J(X)$ containing the polynomial identities of $B_{\infty}$
is generated by a finite number of identities.
\end{theorem}
Later Drensky and Koshlukov \cite{DKo2} found a complete description, with explicitly given polynomial identities,
of the subvarieties of the variety of unitary algebras generated by $B_{\infty}$.

The answer to Problem \ref{problems of Shestakov} (ii) in the case of characteristic 0 follows from the following theorem of Vajs and Zel'manov \cite{VZ}.

\begin{theorem}\label{Specht property for f.g. Jordan algebras}
Let $K$ be a field of characteristic zero. Then every variety generated by a finitely generated Jordan algebra satisfies the Specht property.
\end{theorem}

Medvedev \cite{Me} established an analogue for Jordan algebras of the results of Kruse \cite{Kr} and L'vov \cite{Lv1} for finite associative rings and algebras.

\begin{theorem}\label{PIs of finite Jordan algebras}
Let $\Phi$ be a finite commutative associative ring containing $\frac{1}{2}$.
Then the polynomial identities of every finite Jordan $\Phi$-algebra follow from a finite number.
The variety generated by the algebra has a finite number of subvarieties which are generated by finite algebras.
\end{theorem}

A similar theorem was established by L'vov \cite{Lv3} for finite alternative rings and algebras.

In his Ph.D. Thesis \cite{Is1}, see also \cite{Is2}, Isaev found a basis of the polynomial identities of the algebras $B_d$ and $B_{\infty}$ over a finite field ${\mathbb F}_q$
of characteristic different from 2.

\begin{theorem}\label{Jordan algebra of bilinear form over finite field}
Let the finite field ${\mathbb F}_q$ be of characteristic $p>2$. Then the polynomial identities
\[
([x,y]^2,y,z),\quad ((x-x^q)(y-y^q),z,t),
\]
\[
(x-x^q)(y-y^q)-((x-x^q)(y-y^q))^q,
\]
\[
(x_1,x_2,x_3)(y_1,y_2,y_3)-((x_1,x_2,x_3)(y_1,y_2,y_3))^q
\]
form a basis of the identities of the algebra $B_{\infty}$.
\end{theorem}

Isaev \cite{Is1, Is2} found also the explicit form of the polynomial identities of the algebras $B_d$ over ${\mathbb F}_q$, $\text{char}({\mathbb F}_q)>2$.
It is interesting that for each $d<\infty$ there are two nonisomorphic algebras $B_d$ which have different systems of polynomial identities.

\section{Algebras far from associative}\label{section nonassociative algebras}

Every variety generated by a finite associative, Lie and Jordan algebra satisfies the finite basis property
and has a finite number of subvarieties also generated by finite algebras.
One may expect that the same is true for any finite nonassociative ring or algebra.
But Polin \cite{Po} gave an example of a finite dimensional algebra over a finite field which does not have a finite basis of its polynomial identities.
L'vov \cite{Lv2} gave another simple counterexample:

\begin{theorem}\label{example of Lvov}
Let $K$ be an arbitrary field and let the matrix algebra $M_d(K)$ act canonically from the right on the $d$-dimensional vector space $V_d$.
Consider the algebra $A_d=V_d+M_d(K)$ with multiplication defined by
\[
(v_1+a_1)(v_2+a_2)=v_1a_2,\quad v_1,v_2\in V_d, a_1,a_2\in M_d(K).
\]
{\rm (i)} The algebra $A_2$ does not have a finite basis of its polynomial identities.

{\rm (ii)} The variety $\text{\rm var}(A_3)$ generated by the algebra $A_3$ has an infinite strictly descending chain of subvarieties
whose intersection coincides with the variety $\text{\rm var}(A_2)$.
\end{theorem}

The simplest example in this direction is due to Mal'tsev and Parfenov \cite{MaP} in characteristic 0 and Isaev and Kislitsin \cite{IsKi1}
over an arbitrary field. The approach in \cite{IsKi1} is based on weak polynomial identities.

\begin{theorem}\label{example of Maltsev and Parfenov}
Let $R=V_2+T_2(K)$ be the subalgebra of the algebra $A_2=V_2+M_2(K)$
from Theorem \ref{example of Lvov}, where $T_2(K)$ is the algebra of $2\times 2$ upper triangular matrices.
Then the algebra $R$ does not have a finite basis of its polynomial identities.

{\rm (i)} When the field $K$ is infinite, the basis of the identities of $R$ consists of
\[
x_1(x_2x_3)\text{ and }
x_1[x_2,x_3]y_1\cdots y_n[x_4,x_5],\quad n=0,1,2,\ldots.
\]

{\rm (ii)} Over a finite field $K={\mathbb F}_q$ the basis of the identities of $R$ consists of the identities in {\rm (i)} and the identities
\[
x_1(x_2-x_2^q)(x_3-x_3^q),\quad x_1[x_2,x_3](x_4-x_4^q),\quad x_1(x_2-x_2^q)[x_3,x_4].
\]

{\rm (iii)} Over an arbitrary field $K$ every finite dimensional algebra which contains $R$ as a subalgebra does not have
a finite basis of its polynomial identities.
\end{theorem}

Except the first identity in Theorem \ref{example of Maltsev and Parfenov} the parentheses are left normed, e.g.
\[
x_1x_2^q=((x_1\underbrace{x_2)\cdots)x_2}_{q\text{ times}}\text{ and }x_1[x_2,x_3]=(x_1x_2)x_3-(x_1x_3)x_2.
\]

The subalgebras of the algebra $A_2=V_2+M_2(K)$ have other interesting properties.
Let $V_2$ have a basis $\{v_i\mid i=1,2\}$ and the matrix units $e_{jk}\in M_2(K)$ act on the basis vectors of $V_2$ by the rule
$v_ie_{jk}=\delta_{ij}v_k$, $i,j,k=1,2$, where $\delta_{ij}$ is the Kronecker symbol.
Let $R_1$ and $R_2$ be the subalgebras of $A_2$ with bases $\{v_1,v_2,e_{11},e_{12}\}$ and $\{v_1,v_2+e_{11},e_{21}\}$.
Both algebras satisfy the polynomial identity $x_1(x_2x_3)$.
The following theorem was established by Isaev \cite{Is4}.

\begin{theorem}\label{union of Specht varieties}
Let $K$ be a field of characteristic zero.
Then the varieties of algebras ${\mathfrak R}_1=\text{\rm var}(R_1)$ and ${\mathfrak R}_2=\text{\rm var}(R_2)$
generated, respectively, by the algebras $R_1$ and $R_2$ defined above satisfy the Specht property.
The bases of the polynomial identities of the varieties ${\mathfrak R}_1$ and ${\mathfrak R}_2$ consist of the identity $x_1(x_2x_3)$ and
\[
x_1[x_2,x_3]x_4\text{ for }{\mathfrak R}_1;
\]
\[
x_1[x_2,x_3]x_4-x_4[x_2,x_3]x_1\text{ and }s_3(x_1,x_2,x_3)\text{ for }{\mathfrak R}_2.
\]
The basis of the identities of the union ${\mathfrak R}_1\cup{\mathfrak R}_2=\text{\rm var}(R_1\oplus R_2)$ consists of
\[
x_1(x_2x_3),\quad x_1[x_2,x_3]x_4-x_4[x_2,x_3]x_1,\quad x_1s_3(x_2,x_3,x_4),\quad x_1x_2[x_3,x_4]x_5,
\]
\[
x_1[x_2,x_3]y_1\cdots y_n[x_4,x_5],\quad n=0,1,2,\ldots.
\]
The variety ${\mathfrak R}_1\cup{\mathfrak R}_2$ does not have a finite basis of its identities.
\end{theorem}

Dorofeev \cite{Dor} proved that if two varieties of algebras ${\mathfrak M}_1$ and ${\mathfrak M}_2$
satisfy the descending chain condition for their subvarieties, the same holds for their union ${\mathfrak M}_1\cup{\mathfrak M}_2$.
Here ${\mathfrak M}_1\cup{\mathfrak M}_2$ is the minimal variety containing both ${\mathfrak M}_1$ and ${\mathfrak M}_2$.
Hence the variety ${\mathfrak R}_1\cup{\mathfrak R}_2$ is infinitely based but all of its proper subvarieties have finite bases of identities.

Another result of Isaev and Kislitsin \cite{IsKi2} related with the algebra $A_2=V_2+M_2(K)$
gives an example of a four-dimensional algebra over a finite field without finite basis of its polynomial identities.
It would be interesting to see whether there are three-dimensional algebras with this property.

\begin{theorem}\label{4-dimensional algebra}
Over any finite field ${\mathbb F}_q$ the subalgebra $R$ of $A_2=V_2+M_2({\mathbb F}_q)$ with basis
$\{v_1,v_2,e_{11}+e_{12},e_{22}\}$ does not have a finite basis of its polynomial identities.
The same assertion holds when we consider $R$ as a ring.
\end{theorem}

As in the proof of Theorem \ref{example of Maltsev and Parfenov},
the proofs of Theorems \ref{union of Specht varieties} and \ref{4-dimensional algebra} depend essentially on weak polynomial identities.

We shall mention also the result of Isaev \cite{Is3} that over an arbitrary field $K$
there exists a finite dimensional right alternative algebra without a finite basis of its identities
which answers a question of L'vov \cite[Problem 1.95]{Dn}.
(As we discussed in Section \ref{section Jordan algebras}, by \cite{Lv3} every finite alternative algebra has a finite basis of its identities.)
Finally, see the paper by Oates-MacDonald and Vaughan-Lee \cite{OMVL} for other differences between varieties generated by associative and nonassociative algebras.

\section{L-varieties}\label{section L-varieties}

In this section we shall introduce the notion of L-varieties of pairs and shall survey some results emphasizing on weak polynomial identities of degree 3.

\begin{definition}\label{Omega-ideals}
Let $\Omega\subset K\langle X\rangle$ be a family of polynomials and let
\[
F=\{f_i(x_1,\ldots,x_{n_i})\in K\langle X\rangle\mid i\in I\}.
\]
\begin{itemize}
\item
The pair $(R,V)$ is an $\Omega$-{\it pair} if $\omega(v_1,\ldots,v_n)\in V$ for all $\omega(x_1,\ldots,x_n)\in\Omega$
and all $v_1,\ldots,v_n\in V$.

\item The class of all $\Omega$-pairs satisfying the weak polynomial identities from the system $F$
is the {\it variety of $\Omega$-pairs} defined by the weak polynomial identities from $F$.

\item
The polynomial $f(x_1,\ldots,x_n)\in K\langle X\rangle$ is an $\Omega$-{\it consequence} of the system of polynomials $F$
if $f$ belongs to the ideal of $K\langle X\rangle$ generated by all $f_i(\omega_1,\ldots,\omega_{n_i})$, $\omega_1,\ldots,\omega_{n_i}\in\Omega$.
The system $F$ is a {\it basis of the weak polynomial identities} of the variety of $\Omega$-pairs defined by the identities of $F$.
\end{itemize}
\end{definition}

As we mentioned in the introduction, when $\Omega$ is the free Lie algebra $L(X)$ embedded in $K\langle X\rangle$ and $V$ is a Lie subalgebra of $R^{(-)}$,
we have varieties of Lie pairs. When $\Omega$ is the free special Jordan algebra $SJ(X)$ embedded in $K\langle X\rangle$ and $V$ is a Jordan subalgebra of $R^{(+)}$,
we consider varieties of Jordan pairs. When $\Omega$ is the vector space $KX$ spanned by the free generators $X$ of $K\langle X\rangle$, we have
varieties of arbitrary pairs or L-varieties.
Clearly, in the case $\Omega=K\langle X\rangle$ the pair $(R,R)$ can be identified with the algebra $R$ and we have the usual notion of varieties of associative algebras.

Below, if not explicitly stated, we shall consider the case $\Omega=KX$ only. Then the ideals of weak polynomial identities are invariant under the linear substitutions of the variables,
i.e. under the linear endomorphisms of $K\langle X\rangle$.
{\it If $W$ is such an ideal of $K\langle X\rangle$, it coincides with the ideal of weak polynomial identities of the pair
$(K\langle X\rangle/W,KX)$.} This is an analogue of the well known fact that every ideal $W$ of $K\langle X\rangle$ which is invariant under all endomorphisms of $K\langle X\rangle$
is the T-ideal of the identities of the algebra $K\langle X\rangle/W$.

Our first goal is to survey results on the L-varieties over a field of characteristic 0 satisfying a weak polynomial identity of degree 3.
The following result was stated in 1950 by Malcev \cite{M} where he suggested, independently from Specht \cite{Sp}, to use representation theory of the symmetric group
in the study of PI-algebras. Anan'in and Kemer \cite{AnKe} used this to describe the varieties of associative algebras with distributive lattice of subvarieties.
Since over a field of characteristic 0 every identity is equivalent to a system of multilinear identities, we shall consider the case when $f(x_1,x_2,x_3)$ is multilinear.

\begin{proposition}\label{PIs of degree 3}
Every multilinear polynomial identity of degree $3$ is equivalent to the linearization of one or several of the identities
\begin{equation}\label{identities of degree 3}
x^3,\quad \alpha[x,y]x+\beta x[x,y],\quad \alpha,\beta\in K, (\alpha,\beta)\not=(0,0),\quad s_3(x_1,x_2,x_3).
\end{equation}
\end{proposition}

There are many papers devoted to the description of the varieties of associative algebras satisfying an identity of degree 3:
Anan'in and Kemer \cite{AnKe}, James \cite{Ja}, Klein \cite{Kl}, Nagata \cite{Na}, Olsson and Regev \cite{ORe1, ORe}, Regev \cite{Re1, Re2, Re3}.
The complete description of these varieties was given in the language of lattices of subvarieties by Vladimirova and Drensky \cite{VD}.
Now we shall state the corresponding results for L-varieties comparing them with the results for varieties of (not necessarily unitary) associative algebras.
We shall consider separately each of the identities (\ref{identities of degree 3}).

In 1952 Nagata \cite{Na} proved the following theorem.

\begin{theorem}\label{theorem of Nagata}
Over a field of characteristic $0$ the algebras satisfying the polynomial identity $x^n$ are nilpotent,
i.e. satisfy the polynomial identity $x_1\cdots x_N$ for some $N=N(n)$ depending on $n$.
\end{theorem}

Later Higman \cite{Hi2} generalized this theorem for algebras of characteristic $p>n$ and gave the bound $N(n)\leq 2^n-1$ for the class of nilpotency.
Higman showed also that $N(3)=6$. In the 1980s Gerald Schwarz discovered that Theorem \ref{theorem of Nagata} was established in 1943
by Dubnov and Ivanov \cite{DuIv} but was completely overlooked by the algebraic community.
See the paper by Formanek \cite{F5} for the history of the {\it Dubnov-Ivanov-Nagata-Higman theorem}.
In particular, the value $N(3)=6$ was found by Dubnov \cite{Du} already in 1935.
The better upper and lower bounds for $N(n)$ are due to Razmyslov \cite{Ra4} and Kuz'min \cite{Kuz}:
\[
\frac{1}{2}n(n+1)\leq N(n)\leq n^2.
\]
There is a conjecture that $N(n)=\frac{1}{2}n(n+1)$ and this is confirmed for $n\leq 4$, see the comments in \cite[Part A, Chapter 6]{DF}.
It has turned out that for L-varieties the situation is completely different.
The author of the present paper \cite{D9} constructed an L-variety satisfying the identity $x^3$ which does not have a finite basis of weak identities.

\begin{theorem}\label{WPI xxx}
Over a field of characteristic zero the $L$-ideal generated by the weak polynomial identities
\[
x^3,x_1s_n(x_1,\ldots,x_n)x_2-x_2s_n(x_1,\ldots,x_n)x_1,\quad n\geq 2,
\]
is not finitely generated.
\end{theorem}

The proof uses free products of Grassmann algebras and ideas of the author from 1974
for the construction of nonfinitely based varieties of Lie algebras over a field of positive characteristic \cite{D1}.

Detailed information for the consequences of the polynomial identities of the form
\begin{equation}\label{alpha,beta}
\alpha[x,y]x+\beta x[x,y],\quad \alpha,\beta\in K, (\alpha,\beta)\not=(0,0),
\end{equation}
was obtained by Klein \cite{Kl}, Olsson and Regev \cite{ORe1}, Regev \cite{Re2} and Anan'in and Kemer \cite{AnKe}.
The simplest cases are given in the following theorem of Klein \cite{Kl}.

\begin{theorem}\label{theorem of Klein}
Over a field of characteristic zero every variety of associative algebras which satisfies an identity (\ref{alpha,beta})
satisfies the identity
\[
x_1\cdots x_5-x_{\sigma(1)}\cdots x_{\sigma(5)},\quad \sigma\in S_5,
\]
except the cases when (\ref{alpha,beta}) is equivalent to one of the identities $[x,y]x$, $x[x,y]$ or $[y,x,x]$.
\end{theorem}

The set $\Lambda({\mathfrak M})$ of the subvarieties of any variety $\mathfrak M$ of algebras or of pairs is a lattice with respect to the intersection and union of subvarieties.
By definition, the lattice is {\it distributive} if
\[
{\mathfrak M}_1\cap({\mathfrak M}_2\cup{\mathfrak M}_3)=({\mathfrak M}_1\cap{\mathfrak M}_2)\cup({\mathfrak M}_1\cap{\mathfrak M}_3),
\quad {\mathfrak M}_1,{\mathfrak M}_2,{\mathfrak M}_3\subseteq {\mathfrak M}.
\]
Anan'in and Kemer \cite{AnKe} proved the following theorem.

\begin{theorem}\label{Ananin and Kemer}
The lattice $\Lambda({\mathfrak M})$ of subvarieties of a variety $\mathfrak M$ of associative algebras over a field of characteristic zero
is distributive if and only if $\mathfrak M$ satisfies an identity of the form (\ref{alpha,beta}).
\end{theorem}

The cases when $\alpha+\beta\not=0$ in (\ref{alpha,beta}) are not very complicated and the only difficult case is the identity $[y,x,x]$
which is equivalent to the identity $[x_1,x_2,x_3]$.
It is known that this identity is a basis of the polynomial identities of the Grassmann algebra (Krakowski and A. Regev \cite{KRe}).
See also \cite{VD} where the results are stated in a compact form.

Many important varieties of algebras over a field of characteristic zero have a distributive lattice of subvarieties.
This holds for the variety of unitary associative algebras $\text{var}(M_2(K))$ and the variety of Lie algebras $\text{var}(sl_2(K))$ (Drensky \cite{D10}),
the variety of unitary associative algebras $\text{var}(E\otimes E$) generated by the tensor square of the Grassmann algebra $E$ (Popov \cite{Pop2}),
the variety of unitary Jordan algebras $\text{var}(B_{\infty})$ generated by the algebra of a symmetric bilinear form $B_{\infty}$ (Drensky \cite{D6}).
An analogue of Theorem \ref{Ananin and Kemer} for unitary associative algebras
was established by Popov \cite{Pop1}, with collaboration with Chekova \cite{PopC1, PopC2} and Nikolaev \cite{PopN}.
The result is in the spirit of Theorem \ref{Ananin and Kemer}
but, instead one identity (\ref{alpha,beta}) the variety has to satisfy four identities of degree 5.

The counterpart of Theorem \ref{Ananin and Kemer} for varieties of Lie pairs was established by Drensky and Vladimirova \cite{DV}.

\begin{theorem}\label{distributive Lie pairs}
The lattice $\Lambda({\mathfrak M})$ of subvarieties of a variety $\mathfrak M$ of Lie pairs over a field of characteristic zero
is distributive if and only if $\mathfrak M$ satisfies an identity of the form (\ref{alpha,beta}).
\end{theorem}

As in the case of $x^3$,
the structure of the L-varieties satisfying (\ref{alpha,beta}) is more complicated than the structure of the varieties of associative algebras
satisfying the same identity.

\begin{theorem}\label{WPI (alpha,beta)}
Let $K$ be a field of characteristic zero.

{\rm (i)} Every L-variety satisfying the weak polynomial identity (\ref{alpha,beta})
has a finite basis of its identities.

{\rm (ii)} The lattice of subvarieties of an L-variety is distributive
if and only if it satisfies the weak identity (\ref{alpha,beta}) for some $(\alpha,\beta)\not=(0,0)$.
\end{theorem}

The case $\alpha+\beta=0$ in Theorem \ref{WPI (alpha,beta)} (which is equivalent to the weak polynomial identity $[x_1,x_2,x_3]$)
is due to Volichenko \cite{Vo} as an important step in his proof that the variety ${\mathfrak A}{\mathfrak N}_2$ of Lie algebras
over a field of characteristic zero (which is defined by the identity $[[x_1,x_2,x_3],[x_4,x_5,x_6]]$)
satisfies the Specht property. Volichenko proved even more. He gave a complete description of the weak T-ideals
containing $[x_1,x_2,x_3]$ and closed with respect to affine endomorphisms, i.e. when the set $\Omega$ in the definition of weak T-ideals
is spanned by $K$ and $X$. The case $\alpha=\beta$ which is equivalent to the weak polynomial identity (\ref{WPI for matrices of order 2})
was established by Drensky and Koshlukov \cite{DKo1}, with the complete description of the L-varieties satisfying this weak identity.
As a byproduct \cite{DKo1} contains another proof of Theorem \ref{identities in 2x2 matrices} (i).
Finally, the simpler cases $\alpha\not=\pm\beta$ are handled by the author \cite{D9}.

The weak identities $[x_1,x_2,x_3]$ and (\ref{WPI for matrices of order 2}) have been studied also from other points of view.
Recall that the algebra $M_{p,q}$ consists of block matrices
\begin{equation}\label{element of Mpq}
a=\left(\begin{matrix}a_{11}&a_{12}\\ a_{21}&a_{22}\\ \end{matrix}\right),
\end{equation}
where $a_{11}\in M_p(E_0)$, $a_{12}\in M_{p\times q}(E_1)$, $a_{21}\in M_{q\times p}(E_1)$, $a_{22}\in M_q(E_0)$,
and $E=E_0\oplus E_1$ is the Grassmann algebra with its canonical ${\mathbb Z}_2$-grading,
$M_p(E_0)$ and $M_{p\times q}(E_1)$ are, respectively, the $p\times p$ and $p\times q$ matrices with entries from $E_0$ and $E_1$,
and similarly for $M_q(E_0)$ and $M_{q\times p}(E_1)$.
The algebras $M_{p,q}$ play a key role in the structure theory of T-ideals discovered by Kemer in \cite{Ke5}, see also \cite{Ke3}.
For $M_{p,q}$ the analogue of the usual trace for $M_d(K)$ is the {\it supertrace} defined as
\[
\text{str}(a)=\text{tr}(a_{11})-\text{tr}(a_{22}),
\]
where $a\in M_{p,q}$ is as in (\ref{element of Mpq}).
The following theorem which connects the polynomial identities of $M_{1,1}$ and $E\otimes E$ was established by Kemer \cite{Ke5}
in characteristic 0 and by Azevedo, Fidelis and Koshlukov \cite{AFKo} in positive characteristic.

\begin{theorem}\label{E times E and M11}
{\rm (i)} Over a field of characteristic zero the algebras $M_{1,1}$ and $E\otimes E$ have the same polynomial identities.

{\rm (ii)} Over an infinite field of characteristic $p>2$ the algebras $M_{1,1}$ and $E\otimes E$ have the same multilinear polynomial identities
but $T(M_{1,1})\subsetneqq T(E\otimes E)$. The polynomial identity $[x_1^{p^2},x_2]$ is satisfied by $E\otimes E$
and does not hold for $M_{1,1}$.
\end{theorem}

The weak polynomial identities of $E\otimes E$ were studied by Kemer \cite{Ke4} in his study of nonmatrix polynomial identities.
The description of the weak polynomial identities of $M_{1,1}$ was obtained
by Di Vincenzo and La Scala in \cite{DVLS1} when $\text{char}(K)=0$ and in \cite{DVLS2} in the general case.

\begin{theorem}\label{WPI for M11}
Let $K$ be an infinite field of characteristic different from $2$ and let $W$ be the subspace of $M_{1,1}$ of the elements of supertrace equal to zero.
Then the weak T-ideal of the polynomial identities of the pair $(M_{1,1},W)$ is generated as an ideal invariant under the affine transformations by the weak identities
 \[
[x_1,x_2,x_3]\text{ and } [x_2,x_1][x_3,x_1][x_4,x_1].
\]
\end{theorem}

Vaughan-Lee \cite{VL2} proved that the variety of Lie algebras ${\mathfrak A}{\mathfrak N}_2$ over a field of characteristic 2 does not satisfy the Specht property.
The main step in his proof was the following result.

\begin{theorem}\label{Vaughan-Lee abelian-by-nilpotent}
Over a field of characteristic $2$ the L-variety defined by the weak identities
\[
[x_1,x_2,x_3],\quad [x_1,x_2][x_2,x_3]\cdots[x_{n-1},x_n][x_n,x_1],\quad n=2,3,\ldots,
\]
is not finitely based.
\end{theorem}

It is an open problem whether the variety of Lie algebras ${\mathfrak A}{\mathfrak N}_2$
satisfies the Specht property when $\text{char}(K)=p>2$
and the corresponding problem for the L-variety defined by the weak identity $[x_1,x_2,x_3]$.

The weak polynomial identity $[x_1^2,x_2]$ is involved in the description of the L-variety generated by the pair $(C_d,V_d)$,
where $V_d$ is the $d$-dimensional vector space with a nondegenerate symmetric bilinear form embedded in the Clifford algebra $C_d$.
The following theorem is due to Drensky and Koshlukov \cite{DKo1}.

\begin{theorem}\label{WPI for biliner form}
Let $\text{\rm char}(K)=0$.

{\rm (i)} The weak polynomial identities of the pair $(C_{\infty},V_{\infty})$
follow from the weak identity $[x_1^2,x_2]$ from (\ref{WPI for matrices of order 2}).

{\rm (ii)} The basis of the weak identities of the pair $(C_d,V_d)$, $d<\infty$,
consists of
\[
[x_1^2,x_2]\text{ and } s_{d+1}(x_1,\ldots,x_{d+1}).
\]
\end{theorem}

Now we shall state some results on the L-varieties satisfying the weak standard identity $s_3(x_1,x_2,x_3)$.
Part (i) of the following theorem was established by Kislitsin \cite{Ki1} and part (ii) is due to Isaev and Kislitsin \cite{IsKi1}.
It shows that the L-variety defined by the weak identity $s_3(x_1,x_2,x_3)$ does not satisfy the Specht property over an arbitrary infinite field.
Since the infinite system of weak identities in (ii) consists of multilinear polynomials, it is easy to see that this L-variety does not have a finite basis
also for finite fields.

\begin{theorem}\label{Isaev and Kislitsin s3}
Let the field $K$ be infinite and let ${\mathfrak M}_1$ and ${\mathfrak M}_2$ be the $L$-varieties generated, respectively, by the pairs
$(A_1,A_1)$ and $(A_2,A_2)$, where
\[
A_1=\left(\begin{matrix}
\alpha&\beta\\
0&0\\
\end{matrix}
\right)\text{ and }
A_2=\left(\begin{matrix}
\alpha&0\\
\beta&0\\
\end{matrix}
\right),\quad \alpha,\beta\in K.
\]

{\rm (i)} The $L$-ideals of the weak polynomial identities of ${\mathfrak M}_1$ and ${\mathfrak M}_2$
are generated, respectively, by the weak identities
\[
[x_1,x_2]x_3\text{ and }x_1[x_2,x_3].
\]

{\rm (ii)} The $L$-ideal of the weak polynomial identities of the union
${\mathfrak M}={\mathfrak M}_1\cup{\mathfrak M}_2$ of the L-varieties ${\mathfrak M}_1$ and ${\mathfrak M}_2$
is not finitely generated and is generated by the weak identities
\[
s_3(x_1,x_2,x_3),\quad x_1[x_2,x_3]x_4,\quad [x_1,x_2]y_1\cdots y_n[x_3,x_4],\quad n = 0,1,2,\ldots.
\]
The L-variety $\mathfrak M$ is generated by the pair $(A,V)$, where
\[
A=\left\{\left(\left(\begin{matrix}\alpha&\beta\\ 0&0\\ \end{matrix}\right),\left(\begin{matrix}\alpha&0\\ \gamma&0\\ \end{matrix}\right)\right),\alpha,\beta,\gamma\in K\right\}\subset M_2(K)\oplus M_2(K),
\]
\[
V=\left\{\left(\left(\begin{matrix}\alpha&\beta\\ 0&0\\ \end{matrix}\right),\left(\begin{matrix}\alpha&0\\ \beta&0\\ \end{matrix}\right)\right),\alpha,\beta\in K\right\}.
\]
Hence $\dim(A)=3$ and $\dim(V)=2$.
\end{theorem}

Kislitsin \cite{Ki2} proved the following result.

\begin{theorem}\label{Kislitsin 2017}
{\rm (i)} Over a field of characteristic zero every $L$-ideal containing one of the weak polynomial identities
$[x_1,x_2]x_3$ or $x_1[x_2,x_3]$ is finitely generated.

{\rm (ii)} Over an infinite field of positive characteristic every $L$-ideal $I$, containing one of the weak identities
$[x_1,x_2]x_3$ or $x_1[x_2,x_3]$, is locally finitely generated. This means that for any positive integer $d$
the ideal $I\cap K\langle x_1,\ldots,x_d\rangle$ of $K\langle x_1,\ldots,x_d\rangle$
is a finitely generated as an $L$-ideal.
\end{theorem}

In characteristic 0 the identity $[x_1,x_2]x_3$ is equivalent to the pair of identities $[x_1,x_2]x_1$ and $s_3(x_1,x_2,x_3)$.
The following result of the author \cite{D9} strengthens Theorem \ref{Kislitsin 2017}.

\begin{theorem}\label{Drensky [x,y]z}
{\rm (i)} Over a field of characteristic zero the following $L$-ideals are all $L$-ideals which contain the weak identity
$[x_1,x_2]x_3$:
\begin{itemize}
\item
The $L$-ideal generated by $[x_1,x_2]x_3$;
\item
The $L$-ideal generated by $[x_1,x_2]x_3$ and the weak identity $x_1^n[x_1,x_2]$, $n\geq 0$;
\item
The $L$-ideal generated by $[x_1,x_2]x_3$ and the weak identity $x_1^m$, $m\geq 1$;
\item
The $L$-ideal generated by $[x_1,x_2]x_3$ and the weak identities
$x_1^m$, $m\geq 2$, and $x_1^n[x_1,x_2]$, $0\leq n\leq m-2$.
\end{itemize}
There is also a dual theorem for the weak identity $x_1[x_2,x_3]$.

{\rm (ii)} Over an arbitrary field every $L$-ideal which contains one of the weak polynomial identities
$[x_1,x_2]x_3$ or $x_1[x_2,x_3]$ is finitely generated.
\end{theorem}

The proof of part (i) is based on representation theory of the general linear group
and methods developed by the author in the early 1980s
for the description of the lattices of subvarieties of varieties of linear algebras \cite{D10}.
The proof of part (ii) follows almost immediately from a theorem of Cohen \cite{Coh} from 1967
using the method of Higman \cite{Hi1} from 1952 for partially well ordered sets.
This method is still used to prove the finitely generation of ideals of polynomial identities.

\section{Weak central polynomials}\label{section weak central polynomials}

\begin{definition}\label{weak central polynomials}
The polynomial $c(x_1,\ldots,x_n)\in K\langle X\rangle$ is a {\it weak central polynomial} for the pair $(R,V)$
if $c(x_1,\ldots,x_n)$ is not a weak polynomial identity for the pair and
$c(v_1,\ldots,v_n)$ belongs to the center of $R$ for all $v_1,\ldots,v_n\in V$.
\end{definition}

We shall discuss the following natural problem.

\begin{problem}\label{how to construct weak central}
How to construct weak central polynomials?
\end{problem}

The following theorem is due to Drensky and Zaicev \cite{DZ}.

\begin{theorem}\label{theorem of V.D. and Zaicev}
Let $c(x_1,\ldots,x_n)$ be a multilinear central polynomial for the algebra $R$.
There is an algorithm which, starting with $c(x_1,\ldots,x_n)$, produces a weak central polynomial $c'(x_1,\ldots,x_N)$
for the pair $(R,V)$.
If $\dim(R)=d$ and $\dim(V)=m$, then $\deg(c')\leq n(d-m+1)$.
\end{theorem}

A modification of the algorithm works also for central polynomials of $R$, which are not multilinear.
The most important application of Theorem \ref{theorem of V.D. and Zaicev} concerns
finite dimensional simple (nonassociative) algebras.

We start with a finite dimensional simple (nonassociative) algebra $A$.
We assume that $R=\text{End}_K(A)$ is the associative algebra of the endomorphisms of $A$ as a vector space,
and $V={\mathcal M}(A)$ is the vector space of the operators of left and right multiplication of $A$:
\[
r_a:A\to A, r_a(b)=ba,\ell_a:A\to A, \ell_a(b)=ab,\quad a,b\in A.
\]
By a result of Polikarpov and Shestakov \cite{PSh} $V$ generates the algebra $R$, i.e. $(R,V)$ is a pair.

\begin{remark}\label{disadvantage of the algorithm}
Obviously, every central polynomial of the algebra $R$ is either a weak central polynomial or a weak polynomial identity of the pair $(R,V)$.
The main difficulty in the proof of Theorem \ref{theorem of V.D. and Zaicev} is to find a central polynomial of $R$ which does not vanish evaluated on $V$.
On the other hand
the algorithm in Theorem \ref{theorem of V.D. and Zaicev} has the disadvantage that does not give weak central polynomials
for the pair $(R,V)$ which are not central polynomials for $R$.
For example, let $\dim(R)=d$ and let $R$ have a basis
\[
\{r_1=z_1,\ldots,r_k=z_k,r_{k+1},\ldots,r_d\},
\]
where $\{z_1,\ldots,z_k\}$ is a basis of the center of $R$.
We shall search for a central polynomial of $R$ of the form
\[
c(x_1,\ldots,x_n)=\sum_{\sigma\in S_n}\xi_{\sigma}x_{\sigma(1)}\cdots x_{\sigma(n)}
\]
with $n!$ unknown coefficients $\xi_{\sigma}$.
Then $c(x_1,\ldots,x_n)$ is a cental polynomial if $c(r_{i_1},\ldots,r_{i_n})$ belongs to the center of $R$
for all $d^n$ replacements with elements of the basis of $R$. If
\[
c(r_{i_1},\ldots,r_{i_n})=\sum_{j=1}^d\alpha_{ij}r_j,\quad i=(i_1,\ldots,i_n),
\]
where $\alpha_{ij}$ are linear combinations of $\xi_{\sigma}$, this means that $\xi_{\sigma}$ are solutions of the linear system with $(d-1)d^n$ equations
\[
\alpha_{i1}=\cdots=\alpha_{ik},\alpha_{ij}=0,j=k+1,\ldots,d,
\]
and at least one $\alpha_{i1}$ is different from 0.
If we search for a weak central polynomial for the pair $(R,V)$, $\dim(V)=m$ and $V$ has a basis $\{v_1,\ldots,v_m\}$, then
$c(v_{i_1},\ldots,v_{i_n})$ belongs to the center of $R$ for all $m^n$ replacements with elements of the basis of $V$.
Hence the system which we have to solve has only $(d-1)m^n$ equations.
The naive expectations are that weak central polynomials for the pair $(R,V)$
are more than the central polynomials for $R$.
But the algorithm gives less weak central polynomials.
\end{remark}

\begin{example}\label{weak central polynomials for ad(sl2)}
Let $A=sl_2(K)$, $R=\text{End}_K(sl_2(K))=M_3(K)$, $V=\text{ad}(sl_2(K))$ is the three-dimensional vector subspace of $M_3(K)$ of the multiplications of $sl_2(K)$.
Starting with the central polynomial of Formanek of degree 9, the algorithm gives a weak central polynomial of degree 10.
But the pair
\[
(M_3(K),V)=(\text{End}_K(sl_2(K)),\text{ad}(sl_2(K)))
\]
has weak central polynomials of degree 6:
\[
[x_1,x_2]x_1^3x_2-[x_1,x_2]x_2x_1^3-x_1^3[x_1,x_2]x_2+x_1x_2x_1[x_1,x_2]x_1.
\]
and even a weak central polynomial of degree 3: the standard polynomial $s_3(x_1,x_2,x_3)$.
In the recent paper with Domokos \cite{DoD} we described all weak central polynomials of the pair
\[
(M_3(K),V)=(\text{End}_K(sl_2(K)),\text{ad}(sl_2(K))),\quad \text{char}(K)=0.
\]
\end{example}

\begin{example}\label{weak central for symmetric matrices}
Let $A=H_2(K)$ be the Jordan algebra of the symmetric $2\times 2$ matrices.
Then $V={\mathcal M}(H_2(K))$ is the three-dimensional vector subspace of $M_3(K)$
of the multiplications of $H_2(K)$ and $R=\text{End}_K(H_2(K))=M_3(K)$.
Starting with the central polynomial of Formanek of degree 9
our algorithm gives a weak central polynomial of degree 10.
But the pair $(M_3(K),V)=(\text{End}_K(H_2(K)),{\mathcal M}(H_2(K)))$
has also weak central polynomials of degree eight:
\[
[x_1,x_2,x_1][x_1,x_2][x_1,x_2,x_2]-[x_1,x_2,x_2][x_1,x_2][x_1,x_2,x_1]
\]
\[
+2[x_1,x_2][x_1,x_2][x_1,x_2][x_1,x_2].
\]
\end{example}

Let $U(sl_2({\mathbb C}))$ be the universal enveloping algebra of $sl_2({\mathbb C})$.
Razmyslov \cite{Ra5} studied the weak polynomial identities of the Lie pair $(U(sl_2({\mathbb C})),sl_2({\mathbb C}))$
and of the pairs associated with irreducible representations of $sl_2({\mathbb C})$.

\begin{theorem}\label{irreducible representations of sl2}
{\rm (i)} The variety of Lie pairs generated by the pair $(U(sl_2({\mathbb C})),sl_2({\mathbb C}))$
satisfies the Specht property.

{\rm (ii)} Let $\varrho:sl_2({\mathbb C})\to\text{End}_{\mathbb C}(V_q)\cong M_q({\mathbb C})$
be a $q$-dimensional irreducible representation of $sl_2({\mathbb C})$.
Then the basis of the weak polynomial identities of the Lie pair $(M_q({\mathbb C}),\varrho(sl_2({\mathbb C})))$ consists of three weak polynomial identities:
\[
s_3(x_1,x_2,x_3)x_4-x_4s_3(x_1,x_2,x_3),
\]
\[
\delta\sum_{\sigma\in S_3}\text{\rm sign}(\sigma)[x_4,x_{\sigma(1)},x_{\sigma(2)},x_{\sigma(3)}]-2x_4s_3(x_1,x_2,x_3),
\]
where $\delta=(q^2-1)/4$ is the value of the Casimir element in the representation $\varrho$,
and one more identity in two variables
\[
\text{\rm ART}_q(x_1,x_2)=\text{\rm Ad}(x_2)\prod_{i=1}^{q-1}\left(\ell_{x_2}-\left(i-\frac{1-q}{2}\right)\text{\rm Ad}(x_2)\right)x_1.
\]
Here $\ell_r:R\to R$, $r\in R$, is the operator of left multiplication of the algebra $R$,
and $\text{\rm Ad}(r)(r')=[r,r']$, $r,r'\in R$.
\end{theorem}
For $q=2$ this gives that the weak polynomial identities of the pair $(M_2({\mathbb C}),sl_2({\mathbb C}))$ follow from the weak identity $[x_1^2,x_2]$.

Finally we shall describe some of the results in the Ph.D. thesis of da Silva M\^{a}cedo \cite{dSM} obtained jointly with his advisor Koshlukov, see also \cite{dSMK}.
Recall that the algebra $R$ is a {\it superalgebra} if $R=R_0\oplus R_1$ as a vector space and $R_iR_j\subseteq R_{i+j}$ where $i,j=0,1$ and $i+j$ is taken modulo 2.
Then the {\it Grassmann envelope} of $R$ is the algebra $(R_0\otimes E_0)\oplus(R_1\otimes E_1)$, where $E=E_0\oplus E_1$ is the Grassmann algebra.
The {\it Lie superalgebra} $V=V_0\oplus V_1$ satisfies the {\it super skew-symemtric identity}
\[
[x,y]+(-1)^{\vert x\vert\vert y\vert}[y,x]
\]
and the {\it super Jacobi identity}
\[
(-1)^{\vert x\vert\vert z\vert}[x,[y,z]]+(-1)^{\vert y\vert\vert x\vert}[y,[z,x]]+(-1)^{\vert z\vert\vert y\vert}[z,[x,y]],
\]
where $x,y,z\in V_0\cup V_1$ and $\vert x\vert=0$ or 1 if $x\in V_0$ or $x\in V_1$, respectively.
One of the main results of the theory of Kemer, see \cite{Ke5} and \cite{Ke3}, is that over a field of characteristic 0 any proper variety of associative algebras
is generated by the Grassmann envelope of a finite dimensional associative superalgebra. The following theorem in \cite{dSM, dSMK} gives a partial analogue of the result of Kemer
for varieties of Lie pairs.

\begin{theorem}\label{da Silva Macado and Koshlukov - 1}
Let $\text{\rm char}(K)=0$ and let ${\mathfrak V}=\text{\rm var}(R,V)$ be a variety of Lie pairs such that $R$ is a PI-algebra.
Then there exists a Lie pair $(A,L)=(A_0\oplus A_1,L_0\oplus L_1)$ where both $A$ and $L$ are superalgebras such that $\mathfrak V$ is generated by the pair of Grassmann envelopes
$((A_0\otimes E_0)\oplus (A_1\otimes E_1),(L_0\otimes E_0)\oplus (L_1\otimes E_1))$.
\end{theorem}

As a consequence of Theorem \ref{da Silva Macado and Koshlukov - 1} the authors obtain the following result.

\begin{theorem}\label{da Silva Macado and Koshlukov - 2}
Let the field $K$ be algebraically closed, $\text{\rm char}(K)=0$ and let ${\mathfrak V}=\text{\rm var}(R,V)$ be a variety of Lie pairs such that $R$ is a PI-algebra.
Let $\mathfrak V$ do not contain any pair $(\text{\rm End}(W),\varrho(sl_2(K)))$ for any representation $\varrho:sl_2(K)\to\text{\rm End}(W)$.
Then the Lie algebra $V$ in the pair $(R,V)$ is solvable.
\end{theorem}

We shall finish our paper with the following problem.

\begin{problem}
Describe the weak polynomial identities and the weak central polynomials of the algebra of multiplications of concrete simple nonassociative algebras of small dimension.
\end{problem}

\section*{Acknowledgements}
The author is very grateful to the anonymous referee for the careful reading of the paper and the useful suggestions.

\end{document}